\documentclass[journal]{article}[10pt]
\usepackage{cite}
\usepackage{amsmath,amssymb,amsfonts,amsthm}
\usepackage{hyperref}
\usepackage{algorithmic}
\usepackage{graphicx,graphics}
\usepackage{fancyvrb,placeins,caption}
\usepackage{textcomp}
\usepackage[utf8]{inputenc}
\usepackage{float}
\usepackage{subfigure}

\usepackage{comment}
\usepackage{mathrsfs}
\usepackage{tikz}
\usepackage{url}
\usepackage{color}
\usepackage{lipsum}
\usepackage{geometry}
\usepackage{tabularx,booktabs}
\graphicspath{{pics/}}
\allowdisplaybreaks

\newcolumntype{Y}{>{\centering\arraybackslash}X}

\def\BibTeX{{\rm B\kern-.05em{\sc i\kern-.025em b}\kern-.08em
    T\kern-.1667em\lower.7ex\hbox{E}\kern-.125emX}}

\begin{document}
\title{Electrical Impedance Tomography with Deep Calder\'on Method\thanks{The work of BJ is supported by a start-up fund and Direct Grant for Research 2022/2023, both from The Chinese University of Hong Kong, and that of ZZ by Hong Kong Research Grants Council (15303122) and an internal grant of Hong Kong Polytechnic University (Project ID: P0038888, Work Programme: ZVX3). The work of KS is supported by Basic Science Research Program through the National Research Foundation (Grant No. 2019R1A6A1A11051177) of Korea
(NRF) funded by the Ministry of Education, Korea. The work of ZZ and KS is also supported by CAS AMSS-PolyU Joint Laboratory of Applied Mathematics.}}

\author{Siyu Cen\thanks{Department of Applied Mathematics, The Hong Kong Polytechnic University, Hung Hom, Hong Kong. (\texttt{siyu2021.cen@connect.polyu.hk}, \texttt{zhi.zhou@polyu.edu.hk}).} \and  Bangti Jin\thanks{Department of Mathematics, The Chinese University of Hong Kong, Shatin, N.T., Hong Kong. (\texttt{bangti.jin@gmail.com}, \texttt{b.jin@cuhk.edu.hk}) } \and Kwancheol Shin\thanks{Department of Mathematics, Ewha Womans University, 52, Ewhayeodae-gil, Seodaemun-gu, Seoul
03760, Republic of Korea (\texttt{kcshin3623@gmail.com})} \and Zhi Zhou\footnotemark[2]}

\maketitle

\begin{abstract}
Electrical impedance tomography (EIT) is a noninvasive medical imaging modality utilizing
the current-density/voltage data measured on the surface of the subject. Calder\'on's method
is a relatively recent EIT imaging algorithm that is non-iterative, fast, and capable of reconstructing
complex-valued electric impedances. However, due to the regularization via low-pass filtering and
linearization, the reconstructed images suffer from severe blurring and under-estimation of the exact
conductivity values. In this work, we develop an enhanced version of Calder\'on's method, using {deep} convolution
neural networks (i.e., U-net)  {as an effective targeted post-processing step, and  term the resulting method by deep Calder\'{o}n's method.}
Specifically, we learn a U-net to postprocess the EIT images generated by Calder\'on's method so as to have
better resolutions and more accurate estimates of conductivity values. We simulate chest configurations with which we
generate the current-density/voltage boundary measurements and the corresponding reconstructed images by
Calder\'on's method. With the paired training data, we learn the deep neural network and evaluate its performance
on real tank measurement data. The experimental results indicate that the
proposed approach indeed provides a fast and direct (complex-valued) impedance tomography imaging technique,
and substantially improves the capability of the standard Calder\'on's method.\vskip5pt
\noindent\textbf{Key words}:  Calder\'on's method, electrical impedance tomography, U-net, deep learning.
\end{abstract}

\section{Introduction}\label{sec:introduction}
Electrical Impedance Tomography (EIT) is a noninvasive medical imaging technique utilizing the
electrical property, e.g., conductivity and permittivity, of the concerned subject. The observational
data is measured on the surface of the subject through a number of electrodes attached to the object
by injecting a number of current patterns through the electrodes and measuring the resulting electric voltages at the same time. Based on the measured data, numerical reconstruction algorithms allow
estimating the electric conductivity distribution inside of the target object. This imaging modality enjoys the following
distinct features: radiation-free, low cost and high portability, and thus it has received much attention. The list of successful applications of EIT
includes imaging human thorax, and monitoring the dynamics of lungs and the heart etc. For example, lung tissues
have relatively low conductivity because of the air in lungs, and the heart cells have relatively
high conductivity because of the presence of the muscular tissue and the highly conductive blood inside
the heart. When the subject breaths and the heart beats, the amount of air and blood in each region
changes, and accordingly, the conductivity changes dramatically. Therefore, EIT can be a good imaging
modality for cardio-pulmonary monitoring \cite{Put2019, Bach2018, Vogt2016, Tomicic2019}, and can be
used for setting the pressure of mechanical ventilation \cite{bedside2018}. Compared with computed
tomography and ultrasound imaging which are often used in imaging anatomical structures of the target
area, EIT is regarded as a functional imaging modality because it is also useful to reflect
the functional abnormalities in the target area through impedance changes as well as injuries.

Many reconstruction algorithms have been developed in the literature. Since EIT is highly ill-posed, in
order to get noise-robust reconstructions, each method adapts some regularization techniques. Many
traditional imaging approaches are based on incorporating suitable \textsl{priori} information, typically via suitable regularization 
\cite{DobsonSantosa:1994,AmmariKang:2004,JinKhanMaass:2012,GehreJin:2014,AmmariSeo:2016,AmmariSeo:2017, BerettaSantacesaria:2018}, which in practice require discrete approximations of the forward map and then solving
the resulting finite-dimensional minimization problems using an iterative algorithm, commonly gradient type
methods. In this context, a linearised model is often employed to simplify
the inversion process. 
In contrast, direct algorithms, e.g., factorization method \cite{KirschGrinbert:2008}, D-bar method
\cite{Samuli_2000, Isaacson_texp, Mueller2020}, Calder\'{o}n's method \cite{Bikowski,Peter_2017, Peter_2018},
and direct sampling method \cite{ChowItoZou:2014}, do
not need the forward simulation and thus are much faster than iterative approaches. Calder\'on's method is
a linearlized, fast, direct method, and it can reconstruct the region of interest pointwise
\cite{Calderon, Bikowski, Peter_2017, Shin_2020, Shin_2020_2, Shin_2021}.
It employs special type of complex-valued harmonic functions, i.e., complex geometric optics solutions.
Nonetheless, the quality of the reconstruction by Calder\'on's method has only limited resolution.

Therefore, there is an imperative need to develop enhanced versions of Calder\'{o}n's method. In this work,
we propose a method of post-processing EIT images from Calder\'on's method with U-net, termed as deep Calder\'{o}n method, to exploit the rich prior information contained in the training data while also partially retaining the physical knowledge of the problem. Specifically,
we train a deep neural network (DNN) (i.e., U-net \cite{RonnebergerFischer:2015}) on the pairs of
reconstructions by Calder\'{o}n's method and ground-truth impedances, acting as a post processing
(or denoising) step to enhance the reconstructions by Calder\'{o}n's method. We show that the enhanced
Calder\'on's method can give high quality images while taking advantage of low computational cost
and the ability of handling \textit{complex-valued} impedance. We illustrate the method with four
different examples that simulate a human thorax and targets with complex-valued conductivity that is common in the medical setting. It is
observed that the proposed method not only sharpens the organ boundary in the images but also detects
the correct size and location of the organs. Moreover, the proposed method estimates the conductivity
much more accurately. Also, even very small abnormal inclusions in the lung region can be
detected with the proposed method. Surprisingly, for complex-valued targets, a single trained network
by the real part of the reconstructions can also be used for the post-processing of the  imaginary
part. Since reconstructing one image by Calder\'on's method takes only a few tenth of second, we believe
that it is not only fast to train the network for each case-specific task but also very suitable in real
time imaging once the neural network is trained. In sum, the proposed approach has two salient features. First, partial
physical information is embedded into the initial approximation from Calderon's method, and thus reducing
the computational burden on the neural network (i.e. U-net). Second, the postprocessing via U-net is of
the form image-to-image translation, for which there are many diverse choices, and the employed U-net is quite
effective in the EIT imaging task.

Last we situate the current work in existing literature. Over the last few years, deep learning
techniques have received a lot of attention, and have been
established as the state of art for many medical imaging modalities \cite{Jin2017, Kang2017, Hyun2018}.
This is often attributed to the extraordinary capability of DNNs for approximating high-dimensional
functions, combined with recently architectural and algorithmic innovations and the availability of
a large amount of training data.
The application of DNNs to EIT has also received considerable attention. Much efforts have been put to find direct reconstruction maps from the measured
voltages to conductivity distributions, exploiting the universal approximation properties of DNNs \cite{Hornik}. More recently, more thorough investigations with fully connected
DNN architectures include \cite{Li2019}, and with CNN structures includes \cite{Tan2019, Wu2021}. A
reconstruction algorithm that combines a traditional reconstruction with DNN was implemented in \cite{Martin2017}.
They first reconstructed EIT images by a linearized method, then postprocessed the reconstruction with a trained fully connected DNN, aiming to get noise-robust reconstructed images. A
combination of the D-bar algorithm and U-net network as a post-processing technique is successfully implemented
in \cite{Hamilton2018}, in which the D-bar method is a direct, non-linear reconstruction algorithm.
In contrast, Calder\'on's method is a three-step linearlized D-bar method \cite{Shin_2020}, and therefore much
faster than D-bar method. Also, Calder\'on's method allows reconstructing a complex-valued impedance at
once \cite{Peter_2018, Shin_2020}, without the need of splitting the real and imaginary parts, but its image
resolution compared to the D-bar method is a bit inferior. The work \cite{GuoJiang:2021} proposed a novel
approach to combine the idea of direct sampling method \cite{ChowItoZou:2014} with the neural networks, in
a manner similar to postprocessing, and showed its effectiveness. Fan and Ying \cite{FanYing:2020}
proposed to represent the linearized forward and inverse maps of  EIT with
compact DNNs, whereas Liu et al \cite{WeiLiuChen:2019} proposed an approach based on 
the concept of dominant current.

The rest of this paper is organized as follows. In Section \ref{Sec_elec}, we introduce the mathematical formulation
of the EIT problem, Calder\'on's method and practical implementation. In Section \ref{Deep_cal}, we
develop the enhanced version of Calder\'{o}n's method, termed as deep Calder\'{o}n method, using convolutional neural networks.
The computational results are presented in Section \ref{sec:experiment} with four different examples: a chest phantom with no pathology, chest phantoms with a high or low conductive pathology in one lung, and three slices of cucumber example,
and furthermore we discuss the potential of the proposed method and possible future works. Finally, in Section
\ref{Conclusion}, we give the main conclusions of this work.

\section{Electrical impedance tomography and Calder\'{o}n's method}\label{Sec_elec}
In this section, we give the mathematical formulation of the EIT problem and describe  Calder\'on's method and the computational aspects.
\subsection{The mathematical model for EIT}
Let $\Omega \subset \mathbb{R}^2$ be a bounded domain, $\gamma(x)$ the electrical impedance, and $u(x)$ the electrical potential. Then the governing equation of EIT is given by
\begin{equation}\label{cond_eqn}
\nabla \cdot (\gamma(x) \nabla u(x)) = 0, ~~~ x \in \Omega.
\end{equation}
The current density $j(x)$ and the voltage distribution $f(x)$ on the boundary $\partial\Omega$ are modeled by
\begin{equation*}
\gamma(x)\frac{\partial u}{\partial \nu}(x) = j(x)\quad \mbox{and}\quad u(x)  = f(x),  ~x\in \partial \Omega,
\end{equation*}
respectively, where $\nu$ is the unit outward normal vector to $\partial \Omega$. Therefore, $j(x)$ corresponds to the Nuemann boundary condition and
$f(x)$ corresponds to the Dirichlet boundary condition for problem \eqref{cond_eqn}. The Dirichlet-to-Neumann map (DN map) is then defined as
\begin{equation*}
    \Lambda_{\gamma} : f(x) \longrightarrow \gamma(x) \frac{\partial u}{\partial \nu}(x)|_{\partial\Omega},
\end{equation*}
and the Nuemann-to-Dirichlet map (ND map) is defined as
\begin{equation*}
    \mathcal{R}_{\gamma} : j(x) \longrightarrow  u(x)|_{\partial\Omega}.
\end{equation*}
The forward problem in EIT is to find $\Lambda_{\gamma}$ or $\mathcal{R}_\gamma$ for a given conductivity distribution $\gamma$, i.e.,
$\mathcal{F} : \gamma \rightarrow \Lambda_{\gamma}$, where $\mathcal{F}$ denotes the forward map.
The inverse problem in EIT, also known as Calderon's problem, reads: given the
DN map $\Lambda_\gamma$ or equivalently ND map $\mathcal{R}_\gamma$, find the conductivity 
$\gamma$. That is to invert the map $\mathcal{F}$, i.e., find
the inverse operator $\mathcal{I}=\mathcal{F}^{-1}: \Lambda_{\gamma} \rightarrow  \gamma$. It is well
known that the inverse map $\mathcal{I}$ is highly nonlinear and extremely sensitive to the noise in
the measurement data in the sense that $\mathcal{I}(\Lambda_{\gamma} + \epsilon)$ is very sensitive
to the noise $\epsilon$. Therefore each reconstruction algorithm needs to utilize some regularization
techniques implicitly or explicitly.

\subsection{Calder\'on's method}\label{Calderon_method}
Now we describe Calder\'{o}n's method for EIT reconstruction. The starting point of the method is Calder\'{o}n's pioneering
work \cite{Calderon}. Assuming $\gamma(x) = 1 + \delta(x)$, where $\delta(x) \in L^{\infty}(\Omega)$
is relatively small (i.e., $\gamma(x)$ is a small perturbation from the background conductivity $1$), Alberto Calder\'on \cite{Calderon}
showed that one can recover an approximation to $ \delta(x)$ from the DN map $\Lambda_\gamma$. His short paper showed that this can be done
using Fourier transforms in a clever way, by introducing special solutions to the Laplace equation. Specifically,
Calder\'on utilizes special complex-valued harmonic functions
\begin{align*}
    \phi_1(x;k) = \exp(\pi ik\cdot x + \pi k^{\perp} \cdot x), \quad \mbox{and}\quad
    \phi_2(x;k) = \exp(\pi ik\cdot x - \pi k^{\perp} \cdot x),
\end{align*}
where $\cdot$ denotes Euclidean inner product,  the vector $k = (k_1, k_2)$ and the vector $k^{\perp} = (-k_2, k_1)$, perpendicular to $k$, are
non-physical frequency variables, and $x = (x_1, x_2)$ is the spatial variable for the domain $\Omega$. Note that the magnitudes of
$\phi_1$ and $\phi_2$ grow exponentially to the directions $k^{\perp}$ and $-k^{\perp}$, respectively, and
$$\nabla \phi_1 \cdot \nabla \phi_2 = -2\pi^2|k|^2 \exp(2 \pi i k \cdot x),$$
where the differentiation is with respect to the $x$ variable, and $|\cdot|$ denotes the Euclidean norm of a vector.
For each $i = 1,~2$, and any fixed frequency $k$, consider the conductivity equations for $\omega_i(x) \in H^1(\Omega)$
\begin{align*}
   \left\{\begin{aligned}
   \nabla \cdot (\gamma(x) \nabla \omega_i(x)) &= 0,\quad\mbox{in }\Omega,\\
    \omega_i(x) &= \phi_i(x),\quad \mbox{on }\partial\Omega.
    \end{aligned}\right.
\end{align*}
Since $\gamma(x)$ is a small perturbation from $1$ and $\phi_i(x)$ are harmonic, we can decompose $\omega_i(x)$ as
$$\omega_i(x) = \phi_i(x) + v_i(x),\quad\mbox{with }v_i \in H^1_0(\Omega).$$
Now, we introduce a bilinear form
\begin{equation*}
    B(\phi_1, \phi_2) = \int_{\partial \Omega}\phi_1(x;k) \Lambda_{\gamma} \phi_2(x;k) {\rm d}S.
\end{equation*}
Then, using Green's identity
$$\int_{\Omega}\nabla \phi_1 \cdot \nabla v_2 {\rm d} x= \int_{\partial \Omega}v_2(\nabla \phi_1 \cdot \nu){\rm d}S - \int_{\Omega}v_2 \Delta \phi_1 {\rm d} x =0,$$
we deduce that
\begin{equation}\label{Ek}
\begin{aligned}
     &\int_{\partial \Omega} \phi_1(x;k) \Lambda_{\gamma} \phi_2(x;k) {\rm d}x
    =\int_{\partial \Omega} \omega_1\gamma\nabla \omega_2 \cdot \nu {\rm d}S   \\
    =&\int_{\Omega} \gamma \nabla \omega_1 \cdot \nabla \omega_2 {\rm d }x
    = \int_{\Omega}\gamma\nabla(\phi_1 + v_1) \cdot \nabla(\phi_2 + v_2){\rm d}x  \\
    =&-2\pi^2 |k|^2 \int_{\Omega}\gamma(x)e^{2\pi i x \cdot k }{\rm d} x +\int_{\Omega}\delta(\nabla \phi_1\cdot \nabla v_2 + \nabla \phi_2 \cdot \nabla v_1) + \gamma \nabla v_1  \cdot\nabla v_2 {\rm d}x  \\
    =& -2\pi^2 |k|^2 \int_{\Omega}\gamma(x)e^{2\pi i x \cdot k }{\rm d} x + \Tilde{E}(k)
    = -2 \pi^2 |k|^2 \widehat{\gamma \chi_{\Omega}}(k) + \Tilde{E}(k), 
\end{aligned}
\end{equation}
where $\widehat{\gamma \chi_{\Omega}
}$ is the Fourier transform of $\gamma\chi_{\Omega}$, $\chi_{\Omega}$ is the 
characteristic function of $\Omega$, and the error term $\Tilde{E}(k)$ is given by
\begin{equation*}
    \Tilde{E}(k) = \int_{\Omega}\delta(\nabla \phi_1 \cdot \nabla v_2 + \nabla \phi_2 \cdot \nabla v_1) + \gamma \nabla v_1 \cdot \nabla v_2 {\rm d} x.
\end{equation*}

Therefore, by comparing the first  and last lines of \eqref{Ek}, we get
\begin{align}
\int_{\partial\Omega} \phi_1 \Lambda_{\gamma} \phi_2 {\rm d}x &= -2\pi^2 |k|^2 \widehat{\gamma \chi_{\Omega}}(k) + \Tilde{E}(k)
\quad\mbox{and}\quad
\int_{\partial\Omega}\phi_1 \Lambda_1 \phi_2 {\rm d} x = -2\pi^2 |k|^2 \widehat{\chi_{\Omega}}(k), \label{Lambda1}
\end{align}
where $\Lambda_1$ denote that DN map in the case of $\gamma(x)\equiv 1$, which often 
cannot be measured in practical applications. Upon subtracting each side and noting $\delta=\gamma-1$, we get the Fourier transform {$\widehat{\delta}$} of $\delta$:
\begin{equation*}
\int_{\partial\Omega}\phi_1 (\Lambda_{\gamma} - \Lambda_{1})\phi_2 {\rm d}x = -2 \pi^2 |k|^2 \widehat{\delta}(k) + \Tilde{E}(k).
\end{equation*}
In sum, by dividing both sides by $-2\pi^2|k|^2$ and rearranging terms, we get
\begin{equation*}
    \hat{\delta}(k) = \hat{F}(k) + E(k),
\end{equation*}
where
\begin{align*}
    \hat{F}(k) &= - \frac{1}{2\pi^2|k|^2}\int_{\partial \Omega}\phi_1(\Lambda_{\gamma}-\Lambda_1)\phi_2 {\rm d}x\quad \mbox{and}
   \quad E(k) =\frac{\Tilde{E}(k)}{2\pi^2 |k|^2}.
\end{align*}
Furthermore, Calder\'on \cite{Calderon} showed the following estimate
$$|E(k)| \leq C\|\delta\|^2_{L^\infty(\Omega)}e^{2\pi |k| r},$$
where $C$ is a constant depending on $\Omega$ and $r$ is the radius of smallest sphere containing the domain. It is proven that for a
constant $\alpha$ between $1$ and $2$, if $|k| \leq \frac{2 - \alpha}{2 \pi r} |\log {\|\delta\|_{L^\infty(\Omega)}}|$,
then $|E(k)| \leq C\|\delta\|^{\alpha}_{L^\infty(\Omega)}$ \cite{Bikowski}. Therefore, provided that
$\|\delta\|_{L^\infty(\Omega)}$ and $|k|$ are small enough, $\hat{F}(k)$ is a good approximation to $\hat{\delta}(k)$.
Hence, we may obtain an approximation to $\delta (x)$ by
\begin{equation*}
    \delta(x) \approx \int_{|k| < R} \hat{F}(k)e^{-2 \pi i k \cdot x } {\rm d}k,
\end{equation*}
where $R$ is a constant, controlling the frequency cutoff. Note that truncating $\hat{F}(k)$ to the region
$\{ |k| \leq R\}$ results in a loss of information for high frequency Fourier terms, and hence the obtained
reconstruction tends to be blurry.

Note that $\Lambda_1$ might be approximately measured on some experimental settings such as tank experiments. When we use
measured $\Lambda_1$, we call the reconstructed images \textit{the difference images}. In other practical applications where
$\Lambda_1$ cannot be measured, one needs to simulate $\Lambda_1$ by using numerical techniques, e.g., FEM \cite{Woo1994,GehreJinLu:2014}. Alternatively,
one can utilize (\ref{Lambda1}) and compute the right-hand side integral by using the Simpson's quadrature rule. In both cases, we call
the reconstructed images \textit{the absolute images}. Once $\delta(x)$ is reconstructed, we can get $\gamma = 1+ \delta$ back.

\subsection{Numerical implementation of Calder\'on's method}
In practical applications, we have only a finite number of electrodes attached on the boundary $\partial\Omega$ of a circular domain $\Omega$,  which allow us to collect
partial boundary data. Also, the input current is applied through the electrodes and not current-density continuously along the boundary
and we cannot specify the current densities. Let $L$ be the (even) number of electrodes. Thus we employ the so-called gap model. For $i = 1, 2, \cdots, L-1$, let $\{ T^i\}\subset \mathbb{R}^L$ be a
set of linearly independent current patterns and $\{V^i\}$ the corresponding voltage 
distributions on the electrodes. More precisely, given the $i$th current pattern, let $T^i_l$ and $V^i_{l, \gamma}$, $l = 1, 2,\cdots, L$, denote respectively the applied current and measured voltage on $l$th electrode corresponding the conductivity distribution $\gamma$. By Kirchhoff's law (the law of conservation of charge), we require that $\sum_{l=1}^{L}T^i_l = 0$, and also for a choice of ground, we require $\sum_{l=1}^{L}V^i_{l,\gamma} = 0$. Let $\theta_l= \frac{2\pi 
l}{L}$ be the angle of the midpoint of the $l$th electrode (with respect to the center of the circular domain). In this work, we apply trigonometric current 
patterns defined by
\begin{align}\label{current}
T_{l}^i = \begin{cases}
      M \cos{i \theta_l}, & i = 1, \cdots, \frac{L}{2}-1,\\
      M \cos{\pi l}, & i = \frac{L}{2},\\
      M \sin{(i - L/2)\theta_l)}, & i = \frac{L}{2}+1, \cdots, L-1,
 \end{cases}
\end{align}
where $M$ is the amplitude of the applied current. We denote the normalized current pattern {$t^i_l = T^i_l/{\|T^i\|_{2}}$} and the
corresponding normalized voltage pattern {$v^i_{l, \gamma} = V^i_{l, \gamma}/{\|T^i\|_{2}}$}, where $\|T^i\|_2 = \sqrt{\sum_{l=1}^L (T^i_l)^2}$, and the subscript $\gamma$ indicates the dependence on the conductivity $\gamma$. We model the
current density $j(x)$ on the boundary $\partial\Omega$ by the gab model \cite{Cheney1990}: 
\begin{align*}
    j(x)=\begin{cases}
        \frac{t_l}{A}, & x\in e_l,\\
        0, & \text{otherwise},
    \end{cases}
\end{align*}
where $e_l$ denotes the $l$th electrode and $A$ is the area of the electrode (which is assumed to be of identical area and length for all electrodes).

In order to compute $\hat{F}(k)$ for each $k$, we discretize $\phi_1(x;k)$ and $\phi_2(x;k)$ as column vectors $\phi_1(x_l;k)$ and
$\phi_2(x_l;k)$ and expand them with respect to $\{t^i_l\}$ and $\{ v^i_{l,\gamma}\}$ as
\begin{align*}
    \phi_1(x_l;k) \approx \sum_{i=1}^{L-1} a^i_k t^i_l\quad\mbox{and}\quad
    \phi_2(x_l;k) \approx \sum_{i=1}^{L-1}b^i_k v^i_{l,\gamma}.
\end{align*}
Let $\mathbf{a_k} =(a^1_k,a^2_k, \cdots, a^{L-1}_k)$ and $\mathbf{b_{k,\gamma}} =(b_{k,\gamma}^1, b_{k,\gamma}^2, \cdots, b_{k,\gamma}^{L-1})$ be the vectors with Fourier coefficients with respect to $\{ t^i_l \}_{i=1}^{L-1}$ and $\{v^i_{l,\gamma}\}_{i=1}^{L-1}$, where the subscript $k$ indicates the dependence of the coefficients on the variable $k$. Then,
\begin{align*}
 &\int_{\partial\Omega}\phi_1(x;k) \Lambda_{\gamma}\phi_2(x;k) {\rm d}x
    \approx \int_0^{2\pi} \sum_{i=1}^{L-1} a^i_k t^i(\theta)\left[\Lambda_{\gamma}\sum_{j=1}^{L-1} b^j_{k,\gamma}v^{j,\gamma}(\cdot) \right]{\rm d} \theta\\
    =& \sum_{i=1}^{L-1}\sum_{j=1}^{L-1}a^i_k b^j_{k,\gamma}\int_0^{2\pi}t^i(\theta) \left[ \Lambda_{\gamma}v^{j,\gamma}(\cdot)\right]{\rm d}x
    \approx \frac{1}{A}\sum_{i=1}^{L-1}\sum_{j=1}^{L-1}a^i_k b^j_{k,\gamma} \sum_{l=1}^L t^i_l t^j_l e_w\\
   =& \frac{e_w}{A}\mathbf{a}_k \mathbf{T}\mathbf{b}_{k,\gamma}^T,
\end{align*}
where $e_w$ is the length of the electrodes, $\mathbf{T}$ is a matrix of which $(i,j)$ component is $\sum_{l=1}^L t^i_l t^j_l$, and $\mathbf{b}^T_{k,\gamma}$ is the transpose of $\mathbf{b}_{k,\gamma}$. Similar to the calculation for $\Lambda_1$, one can derive
\begin{equation}
    \hat{F}(k) \approx - \frac{e_w}{2\pi^2A|k|^2}\mathbf{a}_k \mathbf{T}(\mathbf{b}_{k,\gamma} - \mathbf{b}_{k,1})^T.\nonumber
\end{equation}
Alternatively,
\begin{equation}
    \hat{F}(k) \approx - \frac{e_w}{2\pi^2A|k|^2}\mathbf{a}_k \mathbf{T}\mathbf{b}_{k,\gamma}^T - \int_{\Omega}\exp(2\pi i k \cdot x ) {\rm d}x,
\end{equation}
for an absolute image. The inverse Fourier transformation of $\hat{F}(k)$ can be done by FFT. However, we used the Simpson's quadrature rule for the inversion in the experiment.

\section{Deep Calder\'on method}\label{Deep_cal}

 In this section, we propose a
novel approach based on deep neural networks, to enhance the resolution of the image by Calder\'on's method, which is called deep
Calder\'{o}n's method below. It is inspired by the prior works \cite{Jin2017,Hamilton2018,
GuoJiang:2021} which postprocess the initial estimates by direct reconstruction methods (e.g., D-bar method and
direct sampling method) using a trained neural network.

\subsection{Deep Calder\'{o}n method: architecture and training}

In essence, Calder\'on's method reconstructs the target conductivity distribution in the Fourier domain first and then taking
the inverse Fourier transform of the recovered Fourier image $\hat F(k)$. The incompleteness of
the method comes from two distinct sources besides the inevitable measurement and electrode modeling errors. First, in
practical implementation, we only have a finite number of electrodes with which we can only take
partial boundary data. Indeed, if there are $L$ electrodes, the number of independent
measurement and hence the degree of freedom is $L(L-1)/2$ \cite{Cheney1990}. Therefore, from this
limited data, according to Nyquist's sampling theorem, it is impossible to get the infinite precision frequency
data in the frequency domain. Together with the truncation of $\hat{F}(k)$ to the region  $\{|k|<R\}$, this shortage
of information is analogous to the sparse-view reconstruction or under-sampled reconstruction in other imaging modalities (e.g., computed tomoraphy and magnetic resonance imaging)
\cite{Jin2017, Kang2017, Hyun2018}. In these works \cite{Jin2017, Kang2017, Hyun2018}, the problem of under sampling issue is observed to be well overcome by
post-processing the images with a trained CNN. This analogy directly motivates our development of deep Calder\'{o}n method. Second, the other major source of errors stems from neglecting the terms in
\eqref{Ek}, all of which are non-linearly dependent of the unknown inclusion $\delta$, and therefore we only get a linear
approximation $\tilde\delta$ to $\delta$. Since DNNs can approximate any continuous function, a properly trained DNN may extract the
nonlinear relation from the linearly approximated $\tilde\delta$ to the exact $\delta$. Finally, artifacts due to the
measurement noise and electrode modeling errors will be regularized as we trained the network by adding some
amount of noise to the training data \cite{Martin2017}. Consequently, while Calder\'{o}n's method is computationally efficient, it can only provides a rough estimate of
the target image, as most direct reconstruction methods do. Hence, it is of enormous interest to further
improve the resolution of the reconstructions by Calder\'on's method. Following the pioneering works \cite{Jin2017, Kang2017, Hyun2018}, we
propose a deep Calder\'{o}n's method, by postprocessing the image obtained by Calder\'{o}n's method using
deep neural networks that are trained on suitable training data. The experimental evaluation indicates that it is indeed an effective method for EIT reconstruction.

In the proposed deep Calder\'{o}n's method, we look for an image-to-image mapping $f_\theta$ that maps the initial reconstruction $\tilde
\gamma$ by Calder\'{o}n's method to an enhanced image $f_\theta(\tilde \gamma)$ of the same size, where the subscript $\theta$ denotes the collection
of neural network parameters parameterizing the map. This can be viewed as a classical image denoising task
but with specialized noise structure (specific to Calder\'{o}n's method). Following their great successes
in image processing,  we employ convolutional neural networks (CNNs), especially  U-net \cite{RonnebergerFischer:2015}. For the simplicity of discussions,
we assume that the recovered impedance $\tilde \gamma$ is represented as a rectangular image, with $n_x\times n_y$ pixels, where
$n_x$ and $n_y$ denote the numbers of pixels in $x$ and $y$ directions, respectively. Thus the input to
U-net is a matrix $\tilde \gamma \in \mathbb{R}^{n_x\times n_y}$. When the image is of general shape, we
only need to embed it into a rectangular shape. In the experiment, we fill the extended region (outside the domain $\Omega$) with the background conductivity value.

\begin{figure*}[ht]
\includegraphics[width=15cm, height = 7cm]{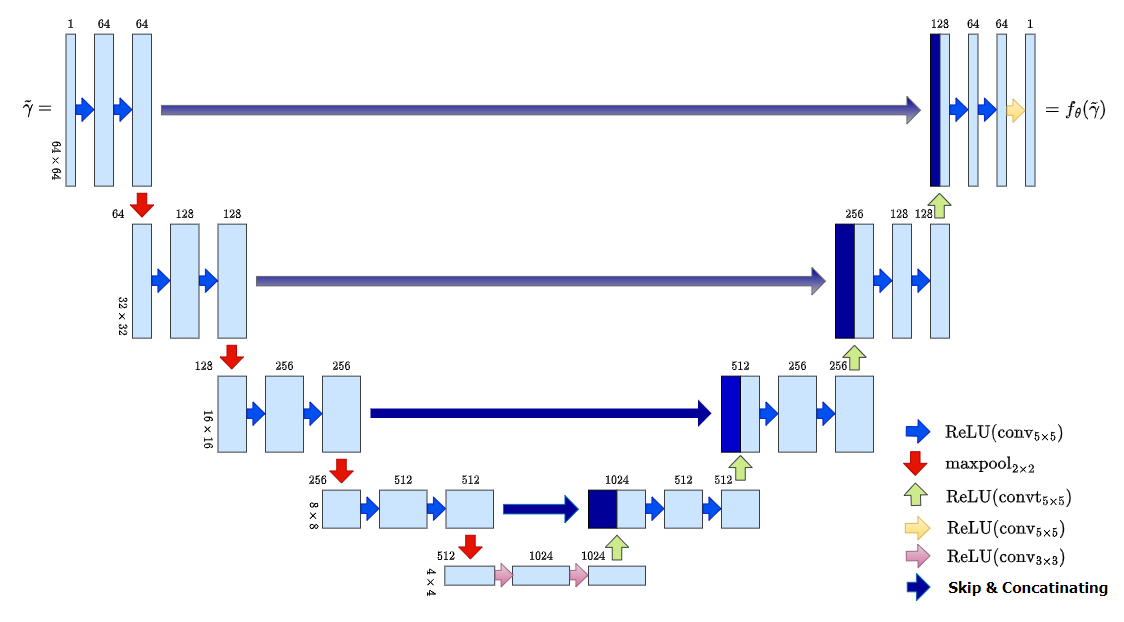}
\caption{The U-net architecture for deep Calder\'on method. The input of the neural network is the reconstruction $\tilde \gamma$ by Calder\'on's method with a resolution $64 \times 64$ and the output is $f_\theta(\Tilde{\gamma})$. Each box represents one layer of the neural network and numbers above the boxes are numbers of channels for each layer. The resolution is on the left bottom of the box. Each convolutional layer is equipped with an activation function ReLU($x$) = max(0,$x$). The blue blocks on the expanding part denote concatenation of blocks from the contracting part (i.e., skip connections). }\label{Unet_diagram}
\end{figure*}

U-net \cite{RonnebergerFischer:2015} is a powerful CNN architecture of encoder-decoder type that consists of a contracting path, a symmetric
expanding path and skip connections; see Fig. \ref{Unet_diagram} for a detailed configuration. We denote the contracting and expanding parts by the superscripts $\rm c$ and $\rm e$, respectively. U-net was
originally developed for image segmentation in biomedical imaging, but with the change of the loss function,
it has also been very successfully applied to regression and image post-processing tasks \cite{Jin2017,Hamilton2018,GuoJiang:2021,Barbano:2022}.

The contracting part consists of several blocks, and each block includes convolution layer, activation
layer and max-pooling layer. The max-pooling layers mainly help extract sharp
features of the input images and reduce the size of the input image by computing the maximum over
each nonoverlapping rectangular region. Specifically, the $i$th block in the contracting part can be
expressed as
\begin{equation*}
  z^{\rm c}_{i+1} := \tau_i^{\rm c} (z_i^{\rm c}) = \mathcal{M}(\rho(W^{\rm c}_i\ast z_i^{\rm c}+b^{\rm c}_i)),
\end{equation*}
where $\ast$ denotes the convolution operation, $W^{\rm c}_i$ and $b^{\rm c}_i$ refer 
to the convolution filter and bias vector, respectively, for the $i$th convolution layer, and 
$z_i^{\rm c}$ is the input image to the $i$th block (in the contracting part, and the input to the first block is the initial reconstruction $\tilde\gamma$ by Calder\'{o}n's method), $\mathcal{M}$
is the max-pooling layer, and $\rho:\mathbb{R}\to\mathbb{R}$ denotes a nonlinear activation function, applied componentwise to a vector. In this study, we employ the rectified linear unit (ReLU) activation function, i.e., $\rho(x)=\max(0,x)$ 
\cite{LeCunBottou:1998}. In practice, one may also include a batch / group
normalization layer that aims to accelerate the training and reduce the sensitivity of the neural network initialization.
Similarly, the expanding part also contains several blocks, each including a transposed 
convolution to extrapolate the input to an image of a larger size. One typical example of the $i$th block is
given by
\begin{equation*}
  z^{\rm e}_{i+1}:=\tau^{\rm e}_i(z_i^{\rm e}) = \mathcal{C}(\rho(\mathcal{T}(z_i^{\rm e},W^{\rm e}_i, b^{\rm e}_i))),
\end{equation*}
where $\mathcal{T}$ denotes the transposed convolution operator, $W^{\rm e}_i$ and $b^{\rm e}_i$ refer
to the corresponding convolutional filter and bias vector at the $i$th layer, and $\mathcal{C}$ is the concatenation
layer. We use $\theta$ to denote the set of all unknown parameters including the convolutional and transposed
convolutional filters and bias vectors, i.e., $\{W_i^{\rm c},b_i^{\rm c}, W_i^{\rm e},b_i^{\rm e}\}$, that have to be learned from the training data via a suitable training procedure.
In addition, we employ a skip connection from the input of the U-net to its output layer at each level, i.e.,
\begin{equation*}
    z_{i+1}^{\rm e} : = \tilde \tau_i^{\rm e}(z_i^{\rm e}) = \tau_{L-i+1}^{\rm c}(z_{L-i+1}^{\rm c}) + \tau_i^{\rm e}(z_i^{\rm e}),
\end{equation*}
where the integer $L$ is the total number of levels in the contracting part of the U-net.
This architectural choice is useful since the input and output images are expected to share similar features, as is in most denoising type tasks. It enforces learning only the
difference between the input and output which avoids learning the (abundantly known) part already contained in the
input. Further, it is known to mitigate the notorious vanishing / exploding gradient problem during the training 
\cite{HeSun:2016}. Then with $L$ contracting blocks
and $L$ expanding blocks, the full CNN model $f_\theta$ can be represented as
\begin{equation*}
  f_\theta(\tilde \gamma) = \tilde \tau_{L}^{\rm e}\circ \cdots \circ \tilde \tau_1^{\rm e} \circ \tau_L^{\rm c}\circ\cdots\circ \tau_1^{\rm c}(\tilde\gamma).
\end{equation*}

To measure the predictive accuracy of the afore-described CNN model $f_\theta$, we employ the mean squared error (MSE) as the loss function
\begin{equation*}
  \mathcal{L}(\theta) = \frac{1}{N}\sum_{n=1}^N \|f_\theta(\tilde \gamma_n) - \gamma_n^{\dag}\|_F^2,
\end{equation*}
where $ \gamma_n^{\dag}$ denotes the true impedance corresponding to the $n$th inclusion sample (of the training dataset, which contains $N$ samples), $\tilde \gamma_n$
denotes the input image also corresponding to the $n$th inclusion sample, obtained by Calder\'{o}n's method, and
$\|\cdot\|_F$ denotes the Frobenius norm of matrices. The
loss $\mathcal{L}(\theta)$ is then minimized by a suitable optimizer, e.g., stochastic gradient descent \cite{RobbinsMonro:1951}, or Adam \cite{KingmaBa:2015} or limited memory BFGS (all of which are directly available in many public software platforms, e.g., PyTorch or TensorFlow), in order to find an (approximate) optimal parameter $\theta^*$. The gradient of the loss $\mathcal{L}(\theta)$ with respect to the U-net parameter $\theta$ can be computed by automatic differentiation.
In this work, the training is carried out in Google Colab engine with the Python library TensorFlow and
the optimization is performed with Adam algorithm with a learning rate $10^{-4}$. The number of trainable parameters is 56,066,369 
when the input image is of size $64\times 64$. The full implementation of the proposed method as well as the data used in the experiments will be made publicly available at the github link \url{https://github.com/KwancheolShin/Deep-Calderon-method.git}.

Conceptually, the proposed deep Calder\'{o}n method can be viewed as a way of incorporating spatial / anatomical \textit{priori} information into the EIT imaging algorithms  \cite{BM1994, Gildewell1995, Vauhkonen1997}. This idea has been explored recently for Calder\'{o}n's method in \cite{Shin_2020, Shin_2021_2}. In the works \cite{Shin_2020, Shin_2021_2}, the approximate location of each organ and their approximate constant conductivity are assumed to be known, since such information can be attained from other imaging modalities, e.g., CT-scans or  ultrasound imaging, and then the prior information is incorporated to the scattering transformation which allows reconstructions with higher-resolution. Note that the method proposed in this work can be viewed as a new way of utilizing the spatial
$\textit{priori}$ information to Calder\'on's method.

\subsection{Constructing training data}

\begin{figure}[ht]
\centering
\begin{tabular}{ccc}
\includegraphics[height = 4cm]{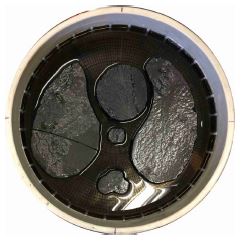}&\includegraphics[height=4cm]{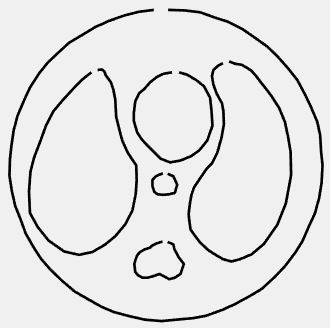}&
\includegraphics[height=4cm]{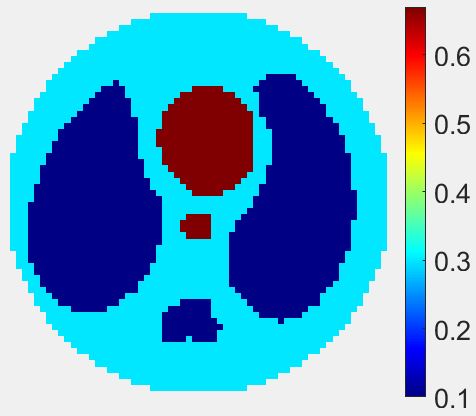}\\
(a) & (b) & (c)
\end{tabular}
\caption{The construction of a realistic conductivity distribution from tank experiment.
(a) The photo of the experimental tank, (b) the boundary of each organ, and (c) the conductivity distribution of the tank from which the EIT data is simulated by FEM solver.}\label{tank}
\end{figure}

The success of the proposed deep Calder\'{o}n's method relies heavily on the availability of paired training data, as most supervised approaches do.
In practice, it is often challenging or expensive to acquire many paired training data experimentally. Hence, we resort to simulated
data. The training of a neural network is done on the simulated data as illustrated in Fig. \ref{tank}. We use Case A in Fig. \ref{4tank} for
the illustration of our method but it is generally all the same for other examples. In order to generate the training
data, from a picture of the tank, we extract the boundary data as in Fig. \ref{tank}(b). Based on the boundary data,
we construct the conductivity distribution $\gamma$ in Fig. \ref{tank}(c). The mesh is on the square $[-1, 1]^2$.

Based on this single reference distribution, we simulate conductivity distributions by varying  the shape and conductivity value.
The detailed method will be described in the experimental study below. We denote these simulated conductivity distributions as
$\{\gamma^\dag_n\}_{n=1}^{N}$, where $N$ is the total number of conductivity distributions. Four instances of the simulated
conductivity distributions $\{\gamma^\dag_n\}_{n=1}^{N}$ for Case A are found in Fig. \ref{4example_No_Rotation}.

With the collection of the simulated conductivity distributions $\{\gamma_n^\dag\}_{n=1}^{N}$, we solve the conductivity equation (\ref{cond_eqn}) with the applied current patterns (\ref{current})
by the FEM \cite{GehreJinLu:2014} to get the simulated current-density/voltage data. Then, we reconstruct an estimate of the conductivity distribution
with Calder\'on's method, as explained in Section \ref{Calderon_method}, on the square $[-1, 1]^2$.  It is well known that EIT is severely ill-posed, and small variations of the conductivity value in the interior of the domain $\Omega$ only lead to very tiny changes in the voltage measurement \cite{Borcea:2002,JinKhanMaass:2012}. In the experiment below, we are interested in very small inhomogeneities located far away from the domain boundary $\partial\Omega$, which can be resolved accurately only if the measurement data is highly accurate, and thus we investigate the addition of $0.01\%$ random
noise to all voltage data below. Accordingly, we set the truncation radius $R$ to $1.3$. We denote the reconstructed images via Calder\'{o}n's method by $\{\tilde\gamma_n\}_{n=1}^{N}$.
In this way, we obtain the paired training data $\{(\tilde\gamma_n, \gamma^\dag_n)\}_{n=1}^{N}$. We randomly split $\{(\tilde \gamma_n,
\gamma^\dag_n)\}_{n=1}^{N}$ into two categories. $90\%$ of the them, denoted by $\{(\tilde\gamma_j, \gamma^\dag_j)\}_{j=1}^{N_1}$,
where $N_1$ is the number of the training data set, are used for the training of the neural network, and the remaining $N_2 = N-N_1$ data
are used for the validation of the neural network.

\section{Experimental evaluation and discussions}\label{sec:experiment}

\begin{figure}[ht]
\centering
\begin{tabular}{cccc}
\includegraphics[height=3.2cm]{ACT4_tank_photo.JPG}&\includegraphics[height=3.2cm]{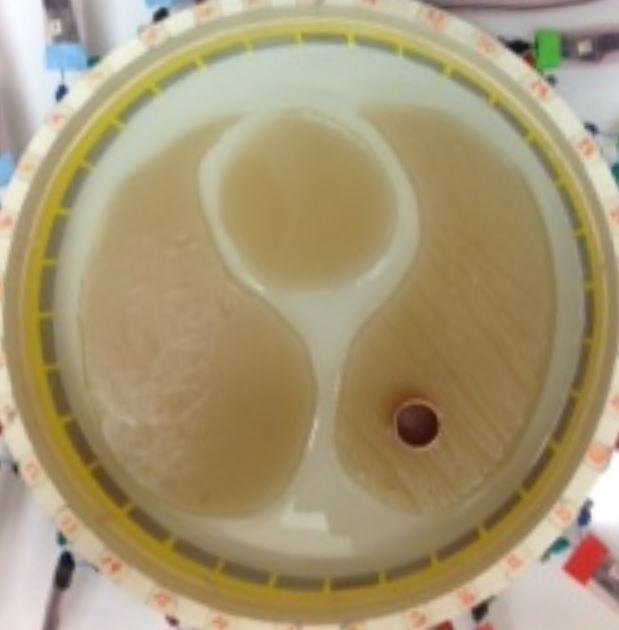}&
\includegraphics[height=3.2cm]{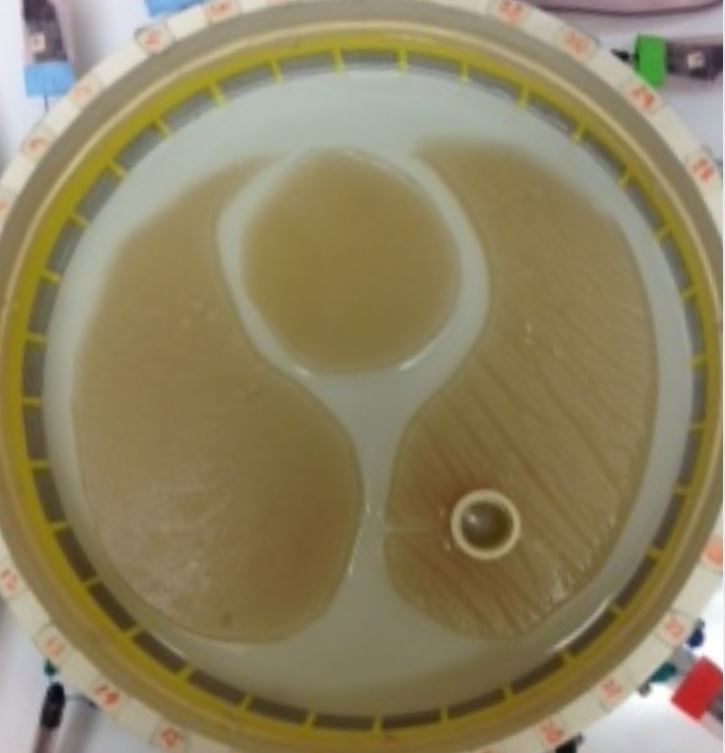}&
\includegraphics[height=3.2cm]{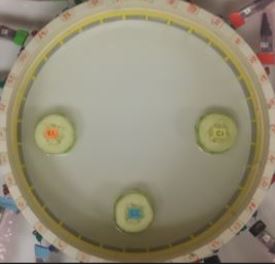}\\
    Case A& Case B & Case C & Case D
\end{tabular}
\caption{The illustration of four experimental settings. From the left to the right, Case A: a chest phantom with no pathology, Case B: a chest phantom with high conductive pathology in a lung, Case C: a chest phantom with low conductive pathology in a lung, and Case D: three slices of cucumber in a tank. Images are credited to \cite{Hamilton2018, Shin_2020}. }\label{4tank}
\end{figure}

The method is tested with four examples; see Fig. \ref{4tank} for the schematic illustration of the phantoms. These setting are commonly use to evaluate EIT in a medical setting. Case A is
a tank experiment with chest phantom filled with a saline bath, and each inclusion is made of agar based targets with
added graphite to simulate the organs in a thorax; two lungs, a heart, a spine and an aorta. In Case B, the targets are
made of agar with added salt and we put a copper pipe in one lung, simulating a high conductive pathology such as a tumor.
In Case C, we put a PVC pipe in the same spot as in Case B to simulate a low conductive pathology such as an air trapping.
In Case D, we have three slices of cucumber to demonstrate the reconstruction of complex-valued impedance. The method of
simulating the training data and the network learning follow closely that in \cite{Hamilton2018}.

\subsection{Case A: Case with no pathology}
{In this first case, we provide a preliminary validation of the proposed postprocessing strategy for Calder\'{o}n's method.}
To generate the simulated training data, we assign values for the conductivity to the regions of background (0.3), left lung (0.1),
right lung (0.1), the heart (0.67), the aorta (0.67), and the spine (0.1). With this choice, we simulate
$N=4096$ conductivity distributions $\{\gamma_n^\dag\}_{n=1}^N$ by varying the size and the conductivity values of two lungs and the heart.
{Specifically, for the size, we give 10$\%$ random variation to the heart and 15$\%$ random variation to each lung as follows. We first compute the center of each organ by averaging each coordinate of the parameterized boundary points. From the center points, we compute the distance to each boundary points, and then multiply a uniformly distributed factor to the distance. The factors are random and independent for each organ.} For the conductivity, we give 30$\%$ variation to each lung, 20$\%$ for the heart, and 10$\%$ for the background. The mesh size is $64 \times 64$. The training takes only a few minutes with a batch size 256 and it is stopped after 12 epoches.

\begin{figure}[ht!]
\centering
\setlength{\tabcolsep}{0pt}
\begin{tabular}{cc}    \includegraphics[width=.48\textwidth, height = 2.8cm ]{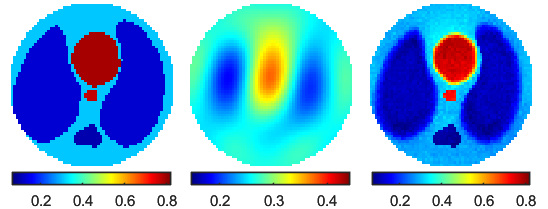}
&\includegraphics[width=.48\textwidth, height = 2.8cm ]{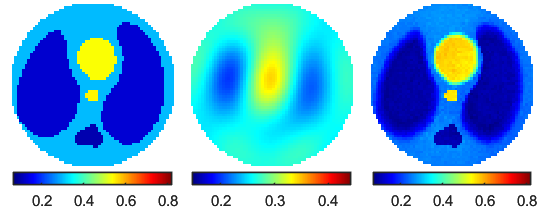}\\
(a) & (b) \\
\includegraphics[width=.48\textwidth, height = 2.8cm ]{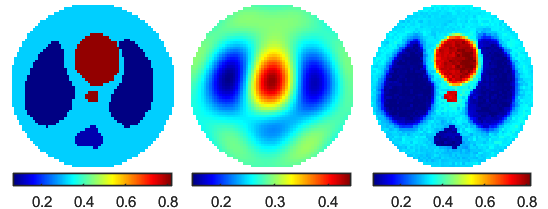}& \includegraphics[width=.48\textwidth, height = 2.8cm ]{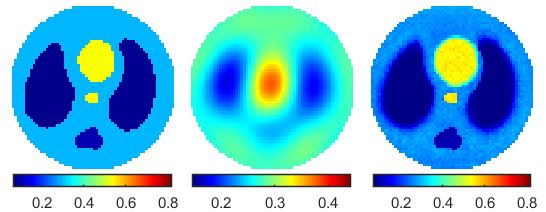}\\
(c) & (d)
\end{tabular}
\caption{Exemplary reconstructions for Case A in four different situations. For each set of images, from left to right, they refer to simulated chest phantom $\gamma_n^\dag$, reconstruction $\tilde\gamma$ by
Calder\'on's method, and the enhanced version $f_\theta(\Tilde{\gamma})$.
}
\label{4example_No_Rotation}
\end{figure}

Fig. \ref{4example_No_Rotation} shows the performance of Deep Calder\'on method on the validation set. Four instances are
displayed to illustrate the results, where in each set of images, we show the simulated phantom $\gamma^\dag$, the
reconstruction $\tilde\gamma$ by Calder\'on's method, and the output $f_{\theta}(\tilde{\gamma})$ from the neural network. Since Calder\'on's method usually does not recover the full
range of the conductivity distribution, the middle column images are under estimated in conductivity values. The four instances show (a) a phantom with big lungs and a high conductive heart, (b) a phantom with big lungs and a relatively low conductive heart, (c) a phantom with
relatively small lungs and a high conductive heart, and (d) a phantom with relatively small lungs and a low conductive
heart. We observe that the trained neural network is able to recover the size of the lungs well, {and the shape of the heart is also well resolved and its size reasonably estimated. More precisely, the ratio of the number of heart cells in $f_{\theta}(\hat{\gamma})$ to that of the ground truth $\gamma^\dag$ is about $0.95, 1.10, 0.93$ and $1.13$ for case (a), (b), (c) and (d), respectively, indicating that the hearts in (a) and (c) are under-estimated in size, whereas the hearts in (b) and (d) are over-estimated in size.} Also, the network is able to recover the full range of conductivity values and also can estimate the approximate conductivity values for both lungs and the heart very well.

\subsection{Cases B and C: Case with a pathology in lungs}
\begin{figure}
\centering
\begin{tabular}{cc}    \includegraphics[width =7cm, height = 3cm ]{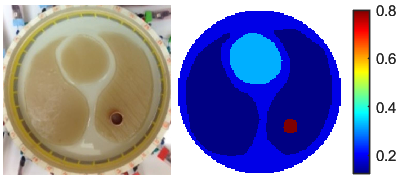} ~~
    &\includegraphics[width =7cm, height = 3cm ]{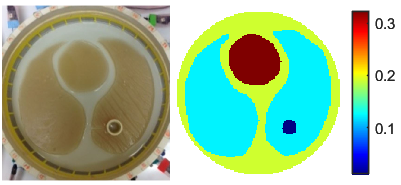}\\
    (a) Case B: phantom with a copper pipe & (b) Case C: phantom with a PVC pipe
\end{tabular}
\caption{Cases B and C: chest phantom with a pipe in one lung and their reference images.}
\label{tankpic}
\end{figure}

{Next we illustrate the approach on the more challenging scenarios, where the lung contains one small inhomeneity, representing pathology or lesion, which is of much interest in practical imaging.}
Cases B and C are simulations with a pathology inside of one lung, cf. Fig. \ref{tankpic}. To simulate the training data,
we assign values for the conductivity to the region of background (0.19), left lung (0.123), right lung (0.123), the heart (0.323),
the copper pipe (0.8), and the PVC pipe (0.01). We regard these distributions, depicted on the right column in Fig. \ref{tankpic},
as the ground truth when we later evaluate the tank experimental data. We simulate $N =16384$ conductivity distributions by randomly
varying the size and the conductivity of each organ. For the size, we give 10$\%$ variation to the heart and 25$\%$ to each lung,
individually. The pathology inclusion is positioned randomly in one lung, {and its size is fixed to be 2.2 cm in diameter throughout}. \textcolor{blue}{In order to locate the inhomogeneity, we first compute  the center $(x_c, y_c)$ of each lung by averaging $x$- and $y$-coordinates of the parameterized boundary points $(x_b, y_b)$, and then use a uniformly distributed random number $r$ between 0 to 1 to locate the center of the inhomogeneity $(x_i, y_i)$ as $(x_i, y_i) = (x_c, y_c) + r\times (x_b - x_c, y_b- y_c)$. In this way, the inhomogeneity may lie inside and on either of the two lungs and can overlap with the lung boundary.} For the conductivity, we give 15$\%$ variation to the heart and 20$\%$ to each lung. One half of the data doesn't include any pathology, 1/4 of the data include a high conductive pathology, and 1/4 of the data include a low conductive inclusion. The mesh size is $128 \times 128$. The training takes about an hour with a batch size $128$ and it is stopped after 14 epoches.

\begin{figure}[ht!]
\centering
\begin{tabular}{cc}
 \includegraphics[width=.48\textwidth, height = 2.75cm ]{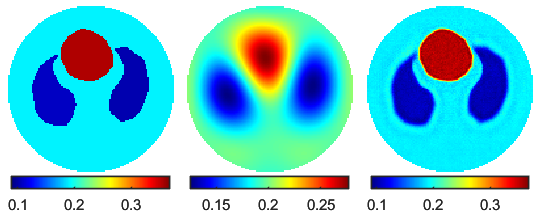}    &  \\
\includegraphics[width=.48\textwidth, height = 2.74cm ]{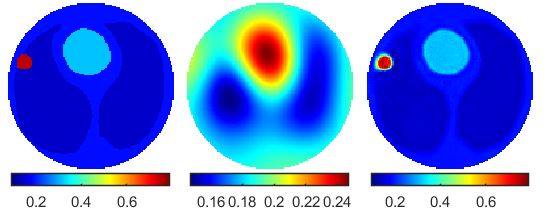}     & \includegraphics[width=.48\textwidth, height = 2.74cm ]{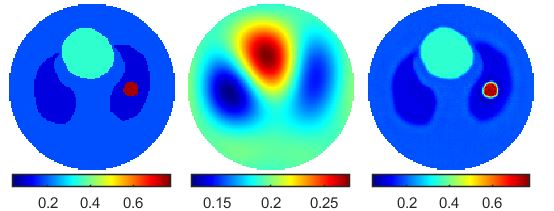}\\
\includegraphics[width=.48\textwidth, height = 2.74cm ]{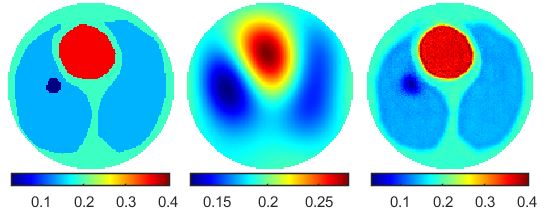}
&    \includegraphics[width=.48\textwidth, height = 2.74cm ]{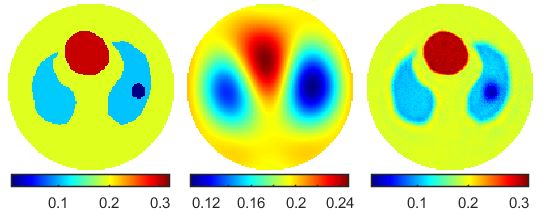}
\end{tabular}
\caption{The reconstruction results with simulated data for Cases B (middle row) and C (bottom row). For each set of images, from left to right, they refer to simulated chest phantom $\gamma^\dag$ (ground truth), the reconstruction $\tilde \gamma$ by Calder\'{o}n's method and enhanced version $f_\theta(\Tilde{\gamma})$.}
\label{copperpipetank_result}
\end{figure}

\begin{figure}
\centering
\begin{tabular}{cc}
 \includegraphics[width=.48\textwidth, height = 2.9cm ]{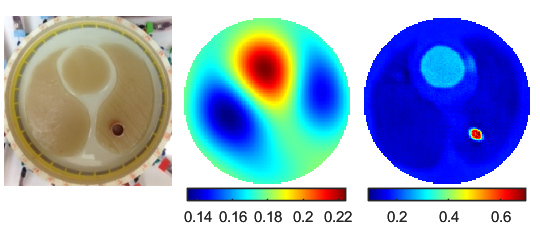}    &     \includegraphics[width=.48\textwidth, height = 2.9cm]{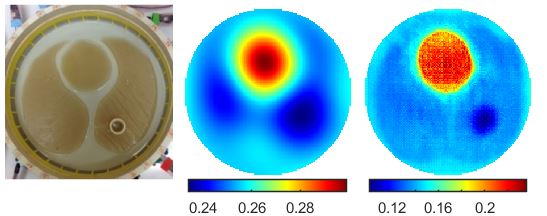} \\
(a) Case B: copper pipe     &  (b) Case C: PVC pipe
\end{tabular}
    \caption{The numerical results (from left to right, ground-truth, reconstruction $\tilde\gamma$ by Calder\'{o}n's method, the enhanced reconstruction $f_\theta(\tilde\gamma)$) with experimental data for Cases B (left) and C (right).}
\label{copperpipetank_result2}
\end{figure}

We present the results on the validation dataset in Fig. \ref{copperpipetank_result}, where the images in the second and third rows are for Cases B and C, respectively and the images in the top row do not contain  small inclusions (pathology). These results show that the network is able to estimate the size and the conductivity of each organ surprisingly well.
Also, the network is able to detect the abnormal inclusion in one lung. Fig. \ref{copperpipetank_result2} shows the result with tank experimental data which were also used in \cite{Alsaker2018, Shin_2020}.
 In each row, the difference images
with Calder\'on's method are shown in the middle column and the corresponding outputs from the network are shown in the last column.
The data were collected using the ACE1 EIT system \cite{ACE1} on a circular tank with diameter 30 cm. The electrodes were 2.54 cm wide,
and the level of saline bath was 1 cm. Adjacent current patterns were applied at the frequency 125 kHz with amplitude 3.3 mA. The copper pipe and
PVC pipe are 1.6 cm and 2.2 cm in diameter. The results in Fig. \ref{copperpipetank_result2} show that the network detects the boundary of each
organ quite well. Also, the copper and PVC pipe is well detected. {Visually, the output $f_{\theta}(\hat{\gamma})$ for Case B by the proposed deep Calder\'{o}n's method appears very fainted} because the range of
reconstruction is very broad due to the presence of the highly conductive copper pipe and this makes the difference in conductivity between
the background and lung region relatively small. In Table 1, from the four reconstructions in Fig. \ref{copperpipetank_result2}, the
conductivity values at four sample positions are displayed. In both cases B and C, the differences in conductivity between inclusions
(H, L, P) and the background (B) are underestimated in $\tilde \gamma$ and are then amplified in $f_\theta(\Tilde{\gamma})$. For
Case B, $\tilde \gamma$, even though the conductivity at (P) is lower than that of (H), for $f_\theta(\Tilde{\gamma})$ the value at (P) is higher
than that of (H), which indicates that the amplification is nonlinear. Given that the amount of variation to the size and the conductivity of inclusions in the training data is quite large and the included
pipes are very small in size for EIT imaging, the results are remarkable.

\begin{table}[ht!]
\centering
\begin{tabular}{| c | c | c | c | c | c | c | }
\cline{4-7}
\multicolumn{3}{c}{} & \multicolumn{2}{|c|}{Case B}
& \multicolumn{2}{c|}{Case C}\\
\hline
 ~region & location & ~conductivity~ & $\tilde \gamma$ & ~$f_\theta(\Tilde{\gamma})$~ & $\tilde \gamma$ & ~$f_\theta(\Tilde{\gamma})$~\\
\hline
 Heart&(41,62) &.323   &.224  & .261  &.199  &.223\\
 Background&(123,65)&.190   &.169  & .167  &.157  &.152\\
 Lung&(71,30) &.123   &.136  & .127  &.141  &.130\\
 Pipe&(89,89) &.8/.01 &.174  & .672  &.134  &.105\\
\hline
\end{tabular}
\caption{Sample conductivity values for Cases B and C from the four reconstructions in Fig. \ref{copperpipetank_result2}. The first column represents the sampling regions: the heart (H), background (B), left lung (L), and the (copper/PVC) pipes (P). The second and third columns show the sampling positions (with respect to (row, column) number from the left-up corner while the size of the images is $124 \times 124$) and the ground truth conductivity values. }\label{table1}
\end{table}

\subsection{Case D: Complex-valued targets}
\begin{figure}[hbt]
\centering
\includegraphics[width=.9\textwidth]{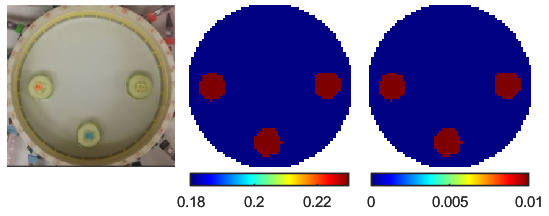}\
    \caption{Case D. The images (from left to right) are the picture of the tank, real and imaginary parts of the approximate ground truth. }
\label{cucumberpic}
\end{figure}

Case D simulates complex-valued targets with three slices of cucumber, see Fig. \ref{cucumberpic} for a schematic illustration of the tank. Each slice of cucumber is about 4.9 cm in
diameter. Due to the cellular structure of cucumber, $\gamma(x) = \sigma(x) +i\omega \epsilon(x)$ is now complex-valued, and we reconstruct both the conductivity $\sigma(x)$ and the permittivity $\epsilon(x)$ simultaneously  when the modulating frequency $\omega$ is known. The conductivity $\sigma(x)$ of the saline bath was 0.180 S/m. We assign values for the complex-valued admitivity to region of cucumbers (0.23 + 0.01$i$) and 0.18 to the region of background.
In the simulation of the training data, three inclusions are randomly positioned (of equal size, 4.9 cm in diameter) without overlapping. {In order to locate the three inclusions, we adopt the polar coordinate system, at the center of the tank, the circular domain $\Omega$ (centered at the origin and normalized to have radius 1). Now for the coordinates $(r, \theta)$ of the center of three inclusions, we randomly pick three radii $r$ each ranging in [0, 0.15], [0.3, 0.45] and [0.6, 0.75] in order to avoid the overlapping of the inclusions but also that they have a good coverage of the domain $\Omega$, and then choose the angle $\theta$ randomly distributed ranging in [0, $2\pi$].} We give 10$\%$ random variation to the conductivity of each inclusion, and does not change the size of the inclusions. The size of the training data set is $N=4096$ on a mesh $64 \times 64$. Training with complex-valued imaging was studied in \cite{Liu}. However, in this study, we train the neural network with only real part, and then use the single trained neural network for the validation of both real and imaginary parts. Also, we normalize the training data to the range $[0, 1]$.

\begin{figure}[ht]
\begin{tabular}{cc}
\includegraphics[width=.5\textwidth, height = 2.75cm]{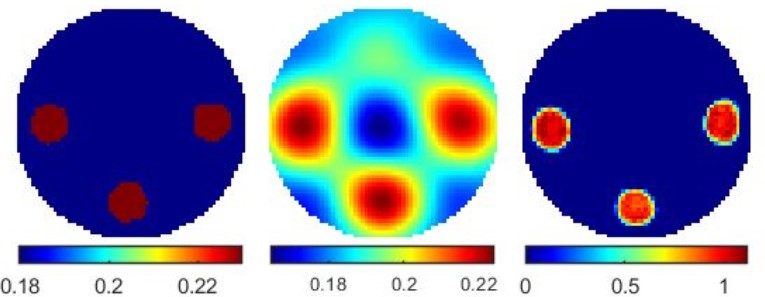} &
\includegraphics[width=.5\textwidth, height = 2.75cm]{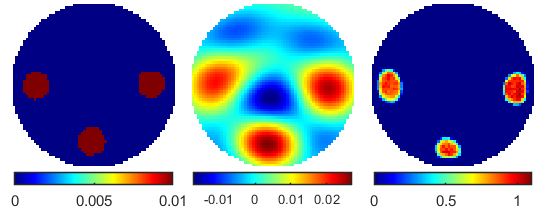}\\
\includegraphics[width=.5\textwidth, height = 3.1cm]{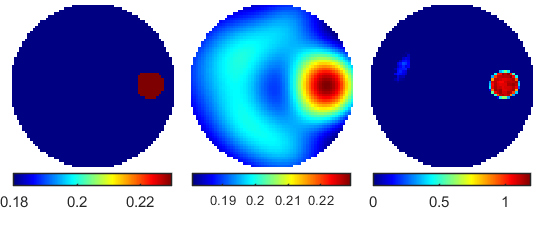} &
\includegraphics[width=.5\textwidth, height = 3.1cm]{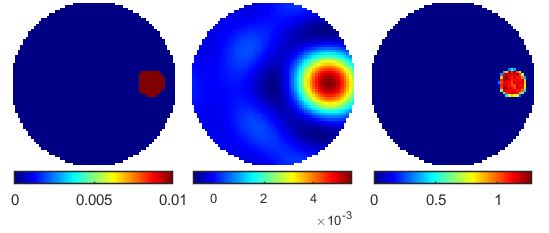}\\
(a) real part & (b) imaginary part
\end{tabular}
\caption{Reconstructions for Case D with experimental tank data (top row), cf. Fig. \ref{cucumberpic} and with simulated data with a single inclusion (bottom row). The images from left to right are the approximated ground truth, Calder\'{o}n reconstruction and the enhanced version.  }\label{tri_ccm_tank}
\end{figure}

We test the proposed method with two examples: tank experimental data and a simulated data. The experimental data was collected using the ACE 1 system at a modulating frequency $\omega=125 \text{ kHz}$, excited with adjacent current
patterns with an amplitude 3.3 mA. The level of saline bath was 1.6 cm. In Fig. \ref{tri_ccm_tank}, we show the approximate ground truth that is extracted from a photo of the tank, reconstruction by Calder\'on's method and the deep version, for the real and imaginary parts separately. Since the training data is normalized to the range $[0, 1]$, when we feed the real part $\Re(\tilde \gamma)$ (or imaginary part $\Im(\tilde \gamma)$) to the trained neural network, we normalize it to $[0,1]$ as well. From Fig. \ref{tri_ccm_tank}, we can see that $f_\theta(\Re(\tilde{\gamma}))$ (or $f_\theta(\Im(\tilde\gamma))$ ranges from 0 to about 1.1 and the location of three inclusions are well detected for both real and imaginary parts.
 For the second test example, we have only one inclusion, cf. the second row of Fig. \ref{tri_ccm_tank}. Even though the case with one inclusion is not in the training dataset, the proposed method can still produce reasonable reconstructions, showing the high out-of-distribution robustness of the proposed approach. This property is highly desirable since in some practical applications, it might be difficult to exactly anticipate the number of inclusions \textsl{a priori} in the training data. 

\subsection{Further discussions}\label{Discussion}
In this work, we have proposed an enhanced Calder\'{o}n method via post-processing using DNNs, called deep Calder\'{o}n's method, and demonstrated its superior performance over classical Calder\'on's method on several medically important settings. Intuitively, this  can be understood as follows. Denote the inversion by Calder\'on's method by
$\mathcal{I}^{\rm cm} : \Lambda_{\gamma} \to \tilde \gamma$,
and the post-processing by U-net by
$f_{\theta} : \tilde \gamma \rightarrow f_\theta(\Tilde{\gamma})$.
While the inverse map $\mathcal{I}: \Lambda_{\gamma} \to \gamma$ is nonlinear in nature, as Calder\'on's method is a linearized and regularized method by the low-pass filtering, the inversion $\mathcal{I}^{\rm cm}$ loses the nonlinear nature of $\mathcal{I}$. Since DNNs are universal approximators \cite{Hornik}, deep Calder\'on's method
$f_{\theta} \circ \mathcal{I}^{\rm cm} : \Lambda_{\gamma} \to \Tilde{\gamma}$
might compensate the loss of non-linearity in $\mathcal{I}^{\rm cm}$ and  therefore it can be a better approximation to $\mathcal{I}$ than the linearised inversion scheme $\mathcal{I}^{\rm cm}$. This effect has been numerically confirmed by the experiments in this work.

The neural network $f_{\theta}$ cannot be independent of the training data and therefore might be different for individual example. Since the structure of a human body is much more complicated than just lungs and the heart, in order to apply the proposed method to a more physically realistic setting including thorax imaging, one must generate a load of training data with arbitrary inclusions to learn the neural network. Moreover, the deep Calder\'{o}n method also depends on many hyper parameters in the experimental setting and Calder\'{o}n's method, e.g., the truncation radius $R$, the injected current pattern and general measurement settings. Hence, constructing a universal neural network is hardly attainable and one will always require some prior knowledge about the target to train a neural network and therefore the suggest method would be case-specific. In this paper, an approximate location, a constant conductivity and the shape of each inclusion are assumed to be known as prior information. The deep Calder\'on method can be viewed as a way of incorporating such implicit prior information into the reconstruction algorithm.
 
In the reconstructions by Calder\'on's method, we have chosen the truncation radius $R$ by a rule of thumb as there is no unified method for choosing $R$ as yet. Having a small $R$ results in a high degree of regularization featuring a stable but blurry reconstruction, and a bigger $R$ leads more detailed but less unstable reconstructed images. Therefore, one may employ a multi-channel network by training the neural network with a large collection of input images for different $R$ simultaneously. The application of deep Calder\'on method to human data involving a moving boundary modeling is also of much interest. 

\section{Conclusions}\label{Conclusion}
We have proposed a deep Calder\'on method and evaluated its performance on several important settings. It is observed that the proposed method can significantly improve the reconstructed images by Calder\'on's method: size, shape and location of the inclusions are well detected, and the boundary of inclusions is sharpened to have resulted in images with much higher resolution. Moreover, the underestimated conductivity values are well adjusted. Also, even a very small abnormal inclusion can be singled out with the proposed method which is almost impossible with the original Calder\'on's method (or other direct reconstrution methods). For complex-valued targets, a single trained neural network using the real part of the reconstructions can be used for the post-processing of both real and imaginary parts. These results indicate that the proposed enhancement strategy is indeed very effective and holds big potentials for certain applications, especially in the medical context.

\bibliographystyle{abbrv}
\bibliography{dcm}

\end{document}